\documentclass{elsarticle}
\usepackage[english]{babel}
\usepackage{a4wide}
\usepackage{mathtools}      
\usepackage{amstext, amsmath, amssymb}
\usepackage{nccmath}
\usepackage{graphicx}
\usepackage{latexsym}
\usepackage{booktabs}
\usepackage{url}
\usepackage{import}
\usepackage{fancyvrb}
\usepackage{stmaryrd}
\usepackage{pdfsync}
\usepackage{algorithm}
\usepackage{tikz}
\usepackage{pgfplots}
\usepackage{pgfplotstable}
\pgfplotsset{compat=newest}
\usepackage{subcaption}
\usepackage[nice]{nicefrac}
\usepackage{tabularx}
\usepackage{todonotes}
\usepackage{enumitem}
\usepackage{xcolor}
\usepackage{rotating}
\biboptions{numbers,square,sort&compress}

\usepackage{array}
\usepackage{verbatim}

\usepackage[plainpages=false]{hyperref}
\hypersetup{
		bookmarksopen=true,
		bookmarksnumbered=true,
		pdfborder={0 0 0},
		colorlinks=true,
		linkcolor=blue,
		urlcolor=blue	
              }

\newcommand{\bfbeta}{\boldsymbol{\bfbeta}}

\DefineVerbatimEnvironment{code}{Verbatim}{frame=single,rulecolor=\color{blue}}

\usepackage[utf8]{inputenc} 
\usepackage[bottom]{footmisc}
\graphicspath{{figures/}} 
\journal{arXiv}

\begin{document}
\begin{frontmatter}
  
  \title{A COMPARISON BETWEEN MESHLESS RADIAL BASIS FUNCTION COLLOCATION METHOD AND FINITE ELEMENT METHOD FOR SOLVING POISSON AND STOKES PROBLEMS}

\author[boun]{İsmet Karakan\corref{cor1}}
\ead{ismet.karakan@boun.edu.tr}

\author[khas]{Dr. Ceren G\"urkan}
\ead{ceren.gurkan@khas.edu.tr}

\author[boun]{Prof. Dr. Cem Avc\i}
\ead{avci@boun.edu.tr}

\cortext[cor1]{Corresponding author}

\address[boun]{Boğaziçi University, Department of Civil Engineering, Bebek 34342, Beşiktaş/İstanbul}
\address[khas]{Kadir Has University, Department of Civil Engineering, Cibali, Kadir Has Cd., 34083 Cibali/Fatih/İstanbul}

\begin{abstract}
Steady and unsteady Poisson and Stokes equations are solved using mesh dependent Finite Element
Method and meshless Radial Basis Function Collocation Method to compare the performances of these
two numerical techniques across several criteria. The accuracy of Radial Basis Function Collocation
Method with multiquadrics is enhanced by implementing a shape parameter optimization algorithm. For
the time-dependent problems, time discretization is conducted using Backward Euler Method. The
performances are assessed over the accuracy, runtime, condition number, and ease of implementation
criteria. Three kinds of errors were calculated; least square error, root mean square error and maximum
relative error. To calculate the least square error while using meshless Radial Basis Function Collocation
Method, a novel technique is implemented. Imaginary numerical solution surfaces are created and then
the volume between those imaginary surfaces and the analytic solution surfaces is calculated, enabling a
fair error calculation. Lastly, all solutions are put together and solution trends are observed over the
number of solution nodes vs. runtime, accuracy vs. runtime, and accuracy vs. the number of nodes. The
assessment indicates the criteria under which Finite Element Method perform better and those when
Radial Basis Function Collocation Method outperforms its mesh dependent counterpart.
\end{abstract}

\begin{keyword}
  Elliptic problems \sep Continuous Galerkin
  \sep Finite Element Method \sep Radial Basis Function Collocation Method \sep Comparison analysis
\end{keyword}

\end{frontmatter}
\section{Introduction}
 
FEM is a numerical method based on a computational mesh, which is basically the combination of the small pieces (elements, and they are finite, hence the name) discretizing the domain of interest. 
The early foundations of the Finite Element Method (FEM) were laid in the 1940s, where McHenry in \cite{McHenry1943} described a method which may fit into almost all types of two-dimensional stress problems involving elastic materials, and it was followed by \cite{courant1943,argyris1954,clough1956}. 
The research conducted since then has led FEM to be a well-established numerical method with a strong theoretical background to be used for the most difficult science and engineering problems. 
Historical development of FEM is briefly given in \cite{gupta96}.  
The meshless methods, on the other hand, are relatively new in the field of computational modeling; Hardy, in 1971, presented a new approach to develop the equations of topography and other irregular surfaces. 
In the meshless approach, unlike mesh dependent methods such as FEM, the domain is decomposed into nodes rather than finite elements or volumes, and the equations of topography are represented by radial basis functions. 
Hardy’s ”Multiquadric (MQ) analysis” aims to find the equation of the topographical surface, which fits all the significant nodes exactly, in a simple manner\cite{Hardymq}. 
Radial Basis Function Collocation Method (RBFCM) is a relatively new development in the area of numerical modeling proposed by Kansa in 1990. He modified MQ scheme to apply it to partial derivative estimates and partial differential equations  \cite{KansaA,KansaB}.
Similar to FEM, RBFCM has been the subject of intense research because of its ease of implementation and low computational cost.
However, there are some disadvantages, mainly due to the poor-condition of the system matrices generated by radial basis functions\cite{Kassabbook}. 
In spite of these disadvantages, technologies dealing with the matrix inversion are getting better day by day. 
Additionally, there are several studies trying to improve the poor-conditioning of radial basis function based methods such as the Localized Collocation Meshless Method\cite{sarler2006,sarler2008,sarler2013} or Radial Basis Function Finite Difference Method\cite{lingdesu}. 
To have a general overview of various meshless methods and their applications, the reader is referred to \cite{Kassabbook}.

Even though the convergence analyses for RBFCM in \cite{madychNelson1990,schaback1993,powell} and for FEM in \cite{zlamal78,babuska81,ainsworth91} are present in the literature, the methods need to be compared on various numerical examples to evaluate their relative performances in practice. 
Several studies have been carried on for both FEM and RBFCM separately; however, comparison studies between these two methods have been limited. 
Li \textit{et al.} compared FEM and RBFCM on 2D steady-state Poisson’s equation over various domains. 
The test function used was only first-order polynomial for the finite element approximations, and the shape parameter, c, was arbitrarily selected\cite{li2003}. 
Gu and Liu worked on 1D and  2D convection-diffusion problems with several Peclet numbers. 
The comparison between the two methods, however, did not include accuracy. 
The conclusion was that the instability issues at high Peclet numbers could easily be overcome by RBFCM with the techniques proposed in the study\cite{gu2006}.
Later in 2012, Golbabai and Raibei studied Stokes eigenproblem in primitive variables of pressure and velocity. A radial basis function formulation for the problem is provided, and the results with FEM variations are compared in the study\cite{golbabai2012}. The comparison was made with respect to relative error only, being similar to various other studies, the model built in this paper was time-independent. 
Ozgener and Tanbay conducted a comparison studying neutron diffusion equations by investigating the effect of shape parameter on the accuracy, convergence rate, and stability of RBFCM. It is concluded in this paper that FEM is more stable than RBFCM due to the local approximation and weak form characteristics of FEM\cite{tanbay2014}.

A review of the previous studies conducted shows that further work will be beneficial to better understand the relative performances of the two methods. 
In the present research, a detailed quantitative comparison between FEM and RBFCM was provided over the criteria such as accuracy, computational time, and condition numbers considering both two-dimensional steady and unsteady Poisson equation and Stokes equation with the the Dirichlet and the the Neumann boundary conditions. 
The standard way of calculating the global stiffness matrix, the Galerkin method of weighted residuals, is used for the finite element analysis. 
Approximations of the first and second-order were considered.  
RBFCM model presented in this work, on the other hand, optimizes the shape parameter for the minimum root mean square error. 
Accuracy assessment was done considering relative, least square error (LSE), and root mean square error (RMSE) for both methods. 
Condition numbers of matrices for both methods are provided quantitatively. 
For computational speed, runtimes have been compared,
and for ease of implementation, a qualitative conclusion has been drawn. 
To fill the gap left, especially by the papers \cite{li2003, golbabai2012}, this study also considers scenarios that include domains where the nodes are randomly distributed for RBFCM approximations, as well as the unsteady models, which are discretized by Backward Euler Method. The effect of mesh refinements on the accuracy, and thus, convergence behaviors, the condition number of system matrices, the computational time needed to run the models (runtime), shape parameter optimization, the effect of higher-order finite element approximations, the effect of thin plate splines (TPS) as well as MQ as radial basis functions for Poisson problems are some of the considerations taken into in this article. This article aims to present strong and weak points of both RBFCM and FEM, and give an a priori idea on better performing one, depending on the desired outcome from the numerical model.

\section{Poisson and Stokes Equations Analysis} 
Poisson equation is a second-order partial differential equation, and although it is the simplest elliptic partial differential equation, it is used to model a wide range of physical phenomena, from gravitational fields to heat transfer and electrical potential. 
Each of these problems has different coefficients in front of the Laplacian operator, defining the problem or material behavior. 
Stokes equation, on the other hand, governs the fully viscous flow when the convective behaviors are neglected. 
The convective terms can be neglected when the velocity of the flow modeled is very small such as groundwater flow and creeping flow. 
Reynold's number, \textit{Re}, is close to zero, and that is why Stokes equation can be considered as a simplified version of Navier-Stokes equations where the nonlinearity arising from convective terms are not present\cite{wathenbook}.

These two equations are analyzed using both FEM and RBFCM for various examples with the Dirichlet and the Neumann type boundary conditions, and the results are compared over various parameters in this study. 
The domains used were chosen to be a unit square in general. 
For all the examples, analytic solutions are set in order to allow for error analyses.  
Two types of finite element functions are considered, namely the first and the second order. For RBFCM, on the other hand, two types of radial basis functions, MQ and TPS, were investigated. 

\subsection{Methodology - FEM}
\label{MFEM}
We start our methodology with the treatment of the simplest partial differential equation, the steady Poisson equation, by using FEM explained here. 
Consider the strong form of steady Poisson equation,
\begin{align}
	\label{eq:Poisson_strongnumformfem}
			\nabla \cdot (-k\nabla u)=f 
\end{align}
This equation can be reduced to a weak form by multiplying it by a test function $v\in P^{(s)}$, $v|_{\Gamma_D}=0$, where $P^{(s)}$ is the space of continuous polynomials of degree less than or equal to s\footnote{Only $s\leq 2$ is used in this study.}, and taking the integral of both sides together with applying integration by parts rule: 
\begin{align}
	\label{eq:ex2weakinnumform}
			\int_\Omega k\nabla u \cdot \nabla v \textnormal{ }d\Omega= \int_\Omega fv \textnormal{ }d\Omega + \int_{\partial\Omega} gv\textnormal{ }dS
\end{align}
where $g=k\nabla u\cdot\vec{n}^{\, }$. The standard finite-element discretization is followed here. The domain is divided into triangular elements, and sets of nodes are defined on each of these elements\footnote{3 nodes exist at each linear triangular element, and 6 nodes exist at each second-order triangular element.}.
The weak form is, then, written for an arbitrary two-dimensional element as follows:
\begin{align}
	\label{eq:weakinnumform}
			\int_{\Omega^e} \left[k\left(\frac{\partial u_h}{\partial x}\frac{\partial v_h}{\partial x}+\frac{\partial u_h}{\partial y}\frac{\partial v_h}{\partial y}\right)\right] \textnormal{ }d\Omega
			= \int_{\Omega^e} fv_h \textnormal{ }d\Omega + \int_{\partial\Omega^e} gv_h\textnormal{ }dS
\end{align}
where $\Omega^e\subseteq\Omega$ consists of the points of that arbitrary triangular element, and $\partial\Omega^e\subseteq\partial\Omega$ consists of the point or points, if exist, which fall on the boundary of the domain.\footnote{The reader may refer to the upper figures in Figure \ref{fig:lshapedssdomain-set-up} for finite element domains where $\partial\Omega=\Gamma_D\cup\Gamma_N$.} 
Assembling the elemental (local) forms in Equation \ref{eq:weakinnumform} for all the elements in the domain and then writing it in the matrix form, we obtain the global stiffness matrix and the load vector as in Equation \ref{eq:KuF}.
\begin{align}
\label{eq:KuF}
Ku=F
\end{align}
where $K$ is the global stiffness matrix, and $F$ is the global load vector.

Next, the time-dependent term for the unsteady version of Equation \ref{eq:Poisson_strongnumformfem} is taken care of. Consider the unsteady Poisson equation,
\begin{align}
\label{eq:unsPoisson}
	\frac{\partial u}{\partial t}+\nabla \cdot (-k\nabla u)=f 
\end{align}
Calculations for the time-dependent term are carried by the mass matrix, $M$. 
The time-dependent term is discretized as follows using Backward Euler Method:
 \begin{align}
	\label{eq:unPoisson_weakinNumformelem3}
			\frac{\partial u}{\partial t}=\frac{ u_h^n-u_h^{n-1}}{\delta t}
\end{align}
The contribution of this new time-dependent term is reflected in the global system through mass matrix, $M$, as follows:
\begin{align}
(M+\delta tK)u^n=\delta tF+Mu^{n-1}
\end{align}
where $\delta t$ is the time step size, and $n$ is the time step.

As we now introduced the FEM discretization for the first PDE this paper studies, next, we move on to the second PDE (\textit{i.e.}, Stokes equation), of which the strong form is given below:
\begin{align}
	\label{stokesStrong}
	-\nabla^2 \textbf{u}+\nabla& p=\textbf{f}\\
	\label{incompCondStrong}
	\nabla \cdot \textbf{u}&=0
\end{align}
where $\textbf{u}=\langle u_x,u_y \rangle$ is the velocity, and $p$ is the pressure. The strong form is reduced to a weak form by multiplying Equation \ref{stokesStrong} with a vector test function, \textbf{v}, and Equation \ref{incompCondStrong} with a scalar test function, $q\in P^{(s)}$ where $q|_{\Gamma_D}=0$. Assembling the local contributions together similar to Poisson case, the weak form can be written as:
\begin{align}
	\label{weakStokes}
	\int_\Omega \nabla \textbf{u} \colon \nabla \textbf{v} \textnormal{ }d\Omega - \int_\Omega p(\nabla \cdot \textbf{v})\textnormal{ }d\Omega &= \int_\Omega \textbf{f} \cdot \textbf{v}\textnormal{ }d\Omega+\int_{\partial \Omega} (\frac{\partial \textbf{u}}{\partial n}-p\vec{n}^{\, })\cdot \textbf{v} dS\\
	\label{weakStokes2}
\int_\Omega (\nabla \cdot \textbf{u}) q\textnormal{ }d\Omega &= 0
\end{align}
The system in Equations \ref{weakStokes},\ref{weakStokes2} is expressed as a system of linear equations as follows:
\begin{equation}
\label{eq:stokesnumformfem1}
	\begin{pmatrix}
	K&0&-L_x\\
	0&K&-L_y\\
	L_x^T&L_y^T&0
	\end{pmatrix}
	\begin{bmatrix}
	u_x\\
	u_y\\
	p
	\end{bmatrix}
	=
	\begin{bmatrix}
	F_x\\
	F_y\\
	0
	\end{bmatrix}
\end{equation}
where $K$ is the stiffness matrix, $F_x$ and $F_y$ are the load vectors with respect to $x$ and $y$ direction, $L_x$ and $L_y$ are the $x$ and $y$ component of the divergence matrix (the details of these matrices are extensively given in \cite{wathenbook}).

Lastly, consider the strong form of the unsteady Stokes equation:
\begin{align}
\label{eq:unsStokes1}
	\frac{\partial \textbf{u}}{\partial t}-\nabla^2 \textbf{u}+&\nabla p=\textbf{f}\\
	\nabla \cdot \textbf{u}&=0
\end{align}
The time-dependent term in Equation \ref{eq:unsStokes1} is written as in Equation \ref{eq:unPoisson_weakinNumformelem3}. 
Finally, the unsteady system of linear equations is written as:
\begin{equation}
\label{eq:unsstokesnumformfem1}
	\begin{pmatrix}
	M+\delta tK&0&-\delta tL_x\\
	0&M+\delta tK&-\delta tL_y\\
	L_x^T&L_y^T&0
	\end{pmatrix}
	\begin{bmatrix}
	u_x\\
	u_y\\
	p
	\end{bmatrix}^n
	=
	\begin{bmatrix}
	\delta tF_x\\
	\delta tF_y\\
	0
	\end{bmatrix}^n
	+
	\begin{pmatrix}
	M&0&0\\
	0&M&0\\
	0&0&0
	\end{pmatrix}
	\begin{bmatrix}
	u_x\\
	u_y\\
	p
	\end{bmatrix}^{n-1}
\end{equation}
\subsection{Methodology - RBFCM}
Here in this section, we explain how to analyze Poisson and Stokes equations with the second numerical technique this paper studies  (\textit{i.e.}, meshless RBFCM). 
RBFCM in this work only refers to Kansa's method as in \cite{KansaA,KansaB}, and it utilizes radial basis functions to approximate partial differential equations.

Among the other types available in the literature, in this study, we only consider MQ and TPS radial basis functions, which are defined as follows:
\begin{itemize}
\item Multiquadrics: $\varphi (r)=(r^2+c^2)^{\beta/2}$, $\beta>0, \beta\in 2\mathbb{N}+1$
\item Thin Plate Splines: $\varphi (r)=r^{\beta}\ln{r}$, $\beta>0, \beta\in 2\mathbb{N}$
\end{itemize}
where $r$ is the Euclidean distance between any two points and $c$ is the shape parameter. In this work, $\beta$ for MQ is chosen to be $1$ only; however, for TPS, $\beta$ varies in each model. 
Hardy in \cite{Hardymq,Hardymq2} used only MQ radial basis functions in his basic scheme. However, the scheme can be written for all types of radial basis functions.

Let $u(\textbf{x})$ be any function, then $u(\textbf{x})$ can be written as a collection of continuously differentiable radial basis functions, $\varphi_j(\textbf{x})$, as such for two dimensions:
\begin{align}
\label{uRBFCM}
u(x,y)=\sum_{j=1}^N\alpha_j\varphi_j(x-x_j,y-y_j)
\end{align}
where $\alpha$ is the vector of the coefficients for each radial basis function $\varphi_j(x,y)$. The aim of the scheme is to solve $\alpha$ first, and then to solve the unknown $u(x,y)$ by using Equation \ref{uRBFCM}. 

In \cite{KansaA,KansaB}, Kansa expanded this basic scheme to approximate partial derivatives and partial differential equations. Partial derivatives of any function $u(x,y)$ can be written as:
\begin{align}
\label{uRABFCM4}
\frac{\partial^{n+m}}{\partial x^n\partial y^m}[u(x,y)]&=\sum_{j=1}^N\alpha_j\frac{\partial^{n+m}}{\partial x^n\partial y^m}[\varphi_j(x-x_j,y-y_j)]
\end{align}
where $n,m\in\mathbb{N}$ and $N$ is the total number of radial basis function centers (RBF centers), (\textit{i.e.}, $(x_j,y_j)$ points). RBF centers are where the radial basis functions are defined.
In addition, another set of points, \textit{collocation nodes}, are defined to carry the information supplied by the differential equation and the boundary conditions\footnote{In all examples provided in this work, these two sets of points are chosen to be exactly identical.}. In short, RBF centers are where RBFCM discretization is done, and collocation nodes are where the numerical solutions are obtained. The idea behind the RBFCM scheme is to assemble radial basis functions in a matrix form and define the problem of interest as a linear system of equations and then, to solve this linear system of equations to find the unknown $\alpha$. Finally, the unknown $u(x,y)$ is approximated by Equation \ref{uRBFCM}. Plugging in the location of collocation nodes in Equation \ref{uRBFCM}, it is expressed as:
\begin{align}
\label{uiRBFCM}
u_i=\varphi_{ij}\alpha_j
\end{align}

Next, consider the strong form of steady Poisson equation, Equation \ref{eq:Poisson_strongnumformfem}, then the partial form and boundary conditions are defined as follows:
\begin{align}
	\label{eq:partialRBFCMnumform1}
			-k(\frac{\partial^2 u}{\partial x^2}+\frac{\partial^2 u}{\partial y^2})=f \textnormal{ in } \Omega \\\label{eq:partialRBFCMnumform2}
			u=u_D \textnormal{ on } \Gamma_D, \\\label{eq:partialRBFCMnumform3}
			\nabla u\cdot {\vec{n}^{\, }}=g \textnormal{ on } \Gamma_N
\end{align}
where f is the source term, and $k$ is a constant material coefficient. Equation \ref{uRABFCM4} is plugged into Equations \ref{eq:partialRBFCMnumform1},\ref{eq:partialRBFCMnumform2},\ref{eq:partialRBFCMnumform3} at the collocation nodes to obtain the following system:
\[
K=
\left(
  \begin{array}{c}
    \varphi_{ij} \\
    \varphi_{ij}^{\vec{n}^{\, }} \\
    -k(\varphi^{xx}+\varphi^{yy})_{ij} \\ 
  \end{array}
\right)
,
F=
\left(
  \begin{array}{c}
    u_{D,i} \\
    g_i \\
    f_i 
  \end{array}
\right)
\begin{array}{l}
\textnormal{on } \Gamma_D,\,\,i=1,2,...,N_{\Gamma_D}\\
\textnormal{on } \Gamma_N,\,\,i=N_{\Gamma_D}+1,N_{\Gamma_D}+2,...,N_{\Gamma_D\cup\Gamma_N}\\
\textnormal{in } \Omega,\,\,i=N_{\Gamma_D\cup\Gamma_N}+1,N_{\Gamma_D\cup\Gamma_N}+2,...,N\\
\end{array}
\] 
where $K$ is the system matrix, and $F$ is the right hand side vector, $N_{\Gamma_D}$ is the total number of nodes on the Dirichlet boundaries, $N_{\Gamma_N}$ is the total number of nodes on the Neumann boundaries, and $j=1,2,...,N$. The Poisson problem defined in Equations \ref{eq:partialRBFCMnumform1},\ref{eq:partialRBFCMnumform2},\ref{eq:partialRBFCMnumform3} is then expressed as a system of linear equations:
\begin{align}
K\alpha=F
\end{align}

For the unsteady problem in Equation \ref{eq:unsPoisson}, a matrix $M$ is defined where M contains the information from the time derivative term, and it is built as follows:
\begin{align}
\label{timeM_RBFCM}
M=
\left(
  \begin{array}{c}
    \varphi_{ij}
  \end{array}
\right)
\begin{array}{c}
\textnormal{in } \Omega\\
\end{array}
\end{align}
Similar to FEM, the time-dependent term is discretized using Backward Euler Method; hence the unsteady Poisson problem is expressed as:
\begin{align}
(M+\delta tK)\alpha^n=\delta tF+M\alpha^{n-1}
\end{align}
where $\delta t$ is the time step size and $n$ is the time step.

Next, the RBFCM formulation for Stokes equation is explained. Consider the strong form of steady Stokes equation, Equations \ref{stokesStrong},\ref{incompCondStrong}, the partial forms of this system of differential equations is written as:
\begin{align}
\label{sstokespartial1}
	-\frac{\partial^2 u_x}{\partial x^2}-\frac{\partial^2 u_x}{\partial y^2}+\frac{\partial p}{\partial x}&=f_x\\
	\label{sstokespartial2}
	-\frac{\partial^2 u_y}{\partial x^2}-\frac{\partial^2 u_y}{\partial y^2}+\frac{\partial p}{\partial y}&=f_y\\
	\label{sstokespartial3}
	\frac{\partial u_x}{\partial x}+\frac{\partial u_y}{\partial y}=0&
\end{align}
To state the system above as a set of linear equations, Equation \ref{uRABFCM4} is substituted in Equations \ref{sstokespartial1}, \ref{sstokespartial2}, \ref{sstokespartial3}. The structure of the linear system generated by this substitution is almost identical to the structure of Equation \ref{eq:stokesnumformfem1}. However, rather than the unknown vector $\langle u_x,u_y,p\rangle$ in Equation \ref{eq:stokesnumformfem1}, the unknown coefficient vector $\langle \alpha_{u_x},\alpha_{u_y},\alpha_{p}\rangle$ is placed. The matrix $K$ and the vectors $F_x$ and $F_y$ are created parallel to the system matrix, $K$, and the right hand side vector, $F$, as explained in this subsection. Besides, $L_x$ and $L_y$ are defined as the first-order partial derivatives of radial basis functions with respect to $x$ and $y$. 
Rather than utilizing efficient iterative methods such as \textit{projection method} \cite{chorin1967} or \textit{pressure correction method} \cite{thomadakis96}, here we chose to use the structure used in \cite{golbabai2012} since its structure is fairly close to the structure in Equation \ref{eq:stokesnumformfem1}.

The matrix $M$ defined in Equation \ref{timeM_RBFCM} contains the information of the time derivative term in unsteady Stokes equations. Similar to steady RBFCM system being equivalent to Equation \ref{eq:stokesnumformfem1}, unsteady RBFCM system is equivalent to Equation \ref{eq:unsstokesnumformfem1} with the exception that the vectors $\langle u_x,u_y,p\rangle^n$ and $\langle u_x,u_y,p\rangle^{n-1}$ in Equation \ref{eq:unsstokesnumformfem1} are being replaced with $\langle \alpha_{u_x},\alpha_{u_y},\alpha_{p}\rangle^{n}$ and $\langle \alpha_{u_x},\alpha_{u_y},\alpha_{p}\rangle^{n-1}$ respectively.

\begin{figure}[htb]
  \begin{center}
    \begin{minipage}[t]{0.47\textwidth}
    \vspace{0pt}
    \includegraphics[width=1.0\textwidth]{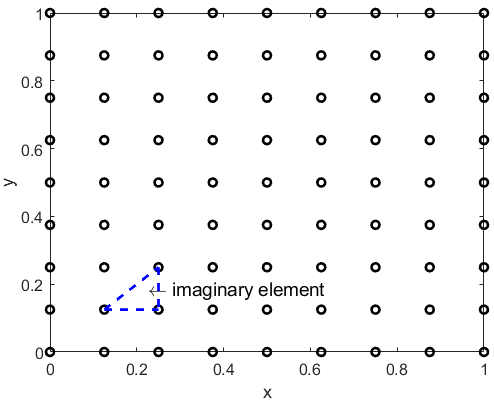}
  \end{minipage}
  \hspace{0.02\textwidth}
  \begin{minipage}[t]{0.47\textwidth}
    \vspace{0pt}
    \includegraphics[width=1.0\textwidth]{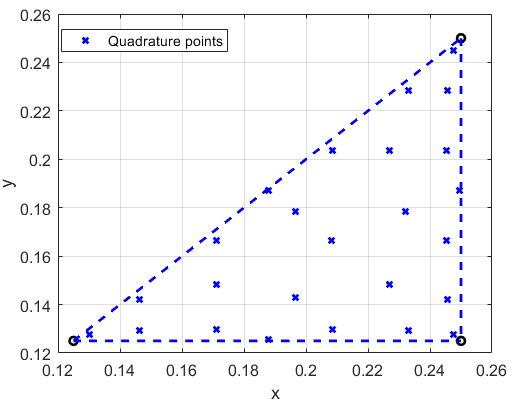}
  \end{minipage}
\end{center}
\caption{An example of an imaginary element in the RBFCM domain (left), Gaussian quadrature points in the imaginary element (right).}
\label{fig:qpsRBFCM}
\end{figure}

Next, the calculation of LSE over the RBFCM setting, (\textit{i.e.}, a setting where there are no "elements" but only "nodes"), is briefly explained. As a novel contribution of this work, we develop an algorithm that mimics the LSE procedure in FEM. 
\textit{Imaginary elements} that correspond to real elements in FEM are generated in the RBFCM domain (see Figure \ref{fig:qpsRBFCM}) to calculate LSE of RBFCM approximations. For each imaginary element, volume integral for the analytical surface and the approximated surface is calculated using Gaussian quadratures. Later, square differences of these volumes are summed for each imaginary element, and the summation is square rooted to calculate LSE.

The volume integral of the approximated surface is calculated by using quadrature points. Assuming a quadrature point has coordinates $(x_p,y_p)$, then unknown at that point, $u_{A}(x_p,y_p)$, is approximated as follows:
\begin{align}
u_{A}(x_p,y_p)=\sum_{j=1}^N\alpha_j\varphi_j(x_p-x_j,y_p-y_j)
\end{align}
Then, the volume integral of the imaginary element is calculated by:
\begin{align}
Volume = \sum_{k=1}^Ew_ku_{Ak}
\end{align}
where $E$ is the total number of quadrature points, and $w_k$ is the Gaussian weights. After finding the volumes, the standard LSE calculation procedure is followed.
\section{Comparison Examples}
\subsection{Steady Poisson Equation with Dirichlet Boundary Conditions}
\label{example1}

For our first numerical example, we have chosen a 2D unit square domain, $\Omega=(0,1)\times(0,1) \subset \mathbb{R}^2$, where all the boundaries are of the Dirichlet type, $\Gamma_D=\partial\Omega$. The analytical solution for this example is:
\begin{align}
\label{eq:ex1analytic}
u(x,y)= \sin(\pi x) \cos(\pi y/2),
\end{align}
To see the behaviors of both methods on a simple trigonometric function, the source term, $f$, then becomes as follows:
\begin{align}
f=\frac{5}{4}\pi^2\sin(\pi x) \cos(\pi y/2),
\end{align}
by setting $k=1$ in Equation \ref{eq:Poisson_strongnumformfem}. To observe the convergence behaviors, the FEM mesh is refined four times, starting from $\delta h=1/4$ down to $\delta h=1/32$.\footnote{Mesh size in x and y direction, $\delta x$ and $\delta y$ is chosen to be equal and referred to as $\delta h$.} Whereas the number of nodes used for RBFCM analysis is set equal to the number of FEM mesh nodes using first-order elements (FEM O(1)). The number of nodes and elements for each $\delta h$ value is presented in Table \ref{dse31}.

\begin{table}[H]
\caption{Domain specifications of Example \ref{example1}.}
\centering
\begin{tabular}{c m{4cm} m{4cm} m{4cm}}
\hline
\boldmath$\delta h$&\textbf{FEM O(1)}&\textbf{FEM O(2)}&\textbf{RBFCM(MQ and TPS)}\\
\hline
1/4&9 nodes in $\Omega$, 16 nodes on $\Gamma_D$, 32 elements.&45 nodes in $\Omega$, 36 nodes on $\Gamma_D$, 32 elements.&9 nodes in $\Omega$, 16 nodes on $\Gamma_D$.\\
1/8&45 nodes in $\Omega$, 36 nodes on $\Gamma_D$, 128 elements.&221 nodes in $\Omega$, 64 nodes on $\Gamma_D$, 128 elements.&45 nodes in $\Omega$, 36 nodes on $\Gamma_D$.\\
1/16&221 nodes in $\Omega$, 64 nodes on $\Gamma_D$, 512 elements.& 961 nodes in $\Omega$, 128 nodes on $\Gamma_D$, 512 elements. &221 nodes in $\Omega$, 64 nodes on $\Gamma_D$.\\
1/32&961 nodes in $\Omega$, 128 nodes on $\Gamma_D$, 2048 elements.& 3969 nodes in $\Omega$, 256 nodes on $\Gamma_D$, 2048 elements. &961 nodes in $\Omega$, 128 nodes on $\Gamma_D$.\\
\hline
\end{tabular}
\label{dse31}
\end{table}

\begin{figure}[H]
  \begin{center}
    \begin{minipage}[t]{0.47\textwidth}
    \vspace{0pt}
    \includegraphics[width=1.0\textwidth]{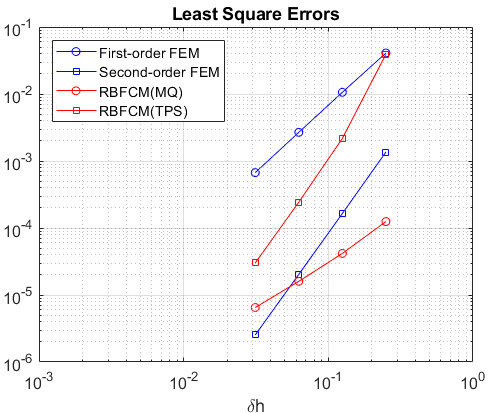}
  \end{minipage}
  \hspace{0.02\textwidth}
  \begin{minipage}[t]{0.47\textwidth}
    \vspace{0pt}
    \includegraphics[width=1.0\textwidth]{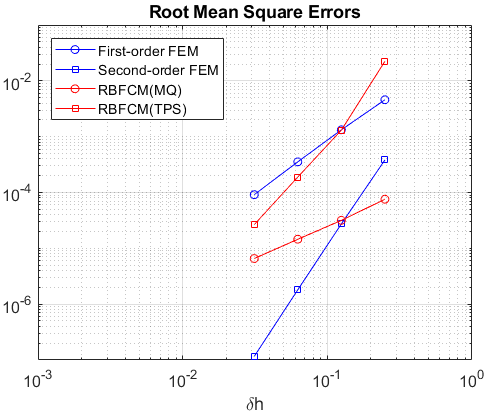}
  \end{minipage}
\end{center}
\caption{Least square errors (left), and root mean square errors (right) of Example \ref{example1}.}
  \label{fig:ex1h-conv}
\end{figure}

LSE and performance results of this example (see Figure \ref{fig:ex1h-conv} and Table \ref{tab:ex1ErrComp}) indicates firstly that both RBFCM-MQ and RBFCM-TPS perform better in accuracy than FEM O(1) for every $\delta h$. The convergence rate (the slope in Figure \ref{fig:ex1h-conv}) of RBFCM-MQ is the least favorable among all methods. It is also interesting to observe that RBFCM-MQ achieved a stable convergence, which is not usual (see \cite{schaback1994}). However, since the shape parameter, $c$, is optimized for RBFCM analysis, it yields the best result at $\delta h=1/4$. The only method that outperforms RBFCM-MQ was the second-order FEM (FEM O(2)) when $\delta h=1/32$. Considering the fact that FEM O(2) utilizes more nodes, it has an excellent convergence rate, and it outperforms, as expected, the other methods for finer meshes. 
On the other hand, the convergence rate of RBFCM-TPS is parallel to that of FEM O(2), even though the gap between the number of nodes that these methods utilize increases as $\delta h$ gets smaller.
Taking a look at the condition numbers of FEM and RBFCM, a clear pattern that RBFCM-MQ always leading to high condition numbers can be observed. It is stated in \cite{schaback1995}, and known as the "Schaback's uncertainty principle", that both condition number of the system matrix and error cannot be small. To acquire more accurate results, either the shape parameter or the number of nodes inside the domain must be increased. However, this principle puts a limit for the error since increasing the shape parameter or the number of nodes leads to higher condition numbers. On the other hand, RBFCM-TPS is a better alternative, with relatively small condition numbers and high accuracy. Finally, FEM has a stronger error convergence behavior with decreases in mesh size, which leads to more accurate solutions.

\begin{table}[H]
\caption[Comparison results for Example \ref{example1}.]{Comparison results for Example \ref{example1}. (LSE: Least Square Error; RMSE: Root Mean Square Error; MRE: Maximum Relative Error; CN: Condition Number; RT: Runtime (sec); OSP: Optimum Shape Parameter.)}
\begin{center}
\begin{tabular}{|c|c|c|c|c|c|c|c|} \hline
	\boldmath$\delta h$&& \textbf{LSE} & \textbf{RMSE}& \textbf{MRE}&\textbf{CN}& \textbf{RT}&\textbf{OSP} \\
	\hline
	1/4&FEM O(1) & 4.091e-02 & 4.578e-03 & 3.099e-02&  1.010e+01 & 0.19&\\
	 &FEM O(2) & 1.361e-03 & 3.826e-04 & 3.175e-03&  3.637e+01 & 5.51&\\
	 &RBFCM-MQ & 1.247e-04 & 7.569e-05 & 3.226e-04&  3.062e+13 & 1.34 & 3.269\\
	 &RBFCM-TPS & 3.922e-02 & 2.239e-02 & 1.187e-01& 5.032e+03 & 0.03 &\\
	 \hline
	 	1/8&FEM O(1) & 1.061e-02 & 1.326e-03 & 1.370e-02&  2.729e+01 &0.78& \\
	 &FEM O(2) & 1.645e-04 & 2.735e-05 & 7.316e-04&1.387e+02 & 22.75&\\
	 &RBFCM-MQ & 4.173e-05 & 3.170e-05 & 5.263e-04&  1.565e+14 & 13.46 & 1.085\\
	 &RBFCM-TPS & 2.172e-03 & 1.309e-03 & 1.688e-02& 6.767e+05 & 0.20 &\\
	 \hline
	1/16&FEM O(1) & 2.678e-03 & 3.542e-04 & 5.098e-03& 1.041e+02 &4.36 &\\
	 &FEM O(2) & 2.034e-05 & 1.809e-06 & 1.785e-04&  5.530e+02 & 94.52&\\
	 &RBFCM-MQ & 1.599e-05 & 1.455e-05 & 6.097e-04&  3.818e+14 & 169.71 & 0.452\\
	 &RBFCM-TPS & 2.430e-04 & 1.892e-04 & 1.045e-02& 5.274e+07 & 2.54 &\\
	\hline
		1/32&FEM O(1) & 6.712e-04 & 9.143e-05 & 1.710e-03& 4.148e+02 &14.21& \\
	 &FEM O(2) & 2.535e-06 & 1.160e-07 & 4.433e-05&  2.213e+03 & 420.99&\\
	 &RBFCM-MQ & 6.510e-06 & 6.564e-06 & 4.427e-03&  2.905e+15 & 2524.41 & 0.218\\
	 &RBFCM-TPS & 3.010e-05 & 2.661e-05 & 9.937e-03& 3.546e+09 & 36.91 &\\
	 \hline
	\end{tabular}
\label{tab:ex1ErrComp}
\end{center}
\end{table}

\subsection{Steady Poisson Equation with Dirichlet and Neumann Boundary Conditions}
\label{lshapedss}

\begin{figure}[H]
  \begin{center}
    \begin{minipage}[t]{0.47\textwidth}
    \vspace{0pt}
    \includegraphics[width=1.0\textwidth]{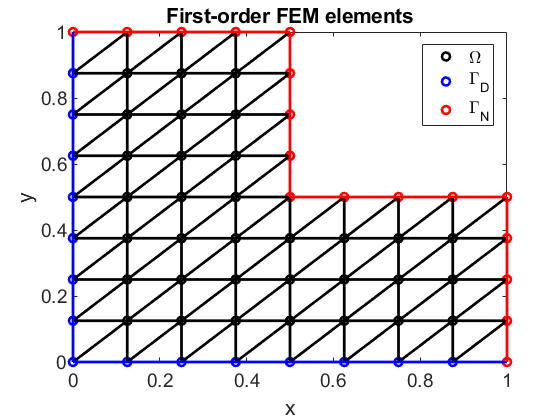}
  \end{minipage}
  \begin{minipage}[t]{0.47\textwidth}
    \vspace{0pt}
    \includegraphics[width=1.0\textwidth]{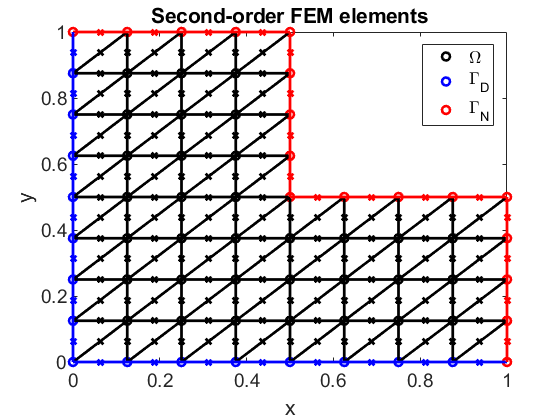}
  \end{minipage}
    \begin{minipage}[t]{0.47\textwidth}
    \vspace{0pt}
    \includegraphics[width=1.0\textwidth]{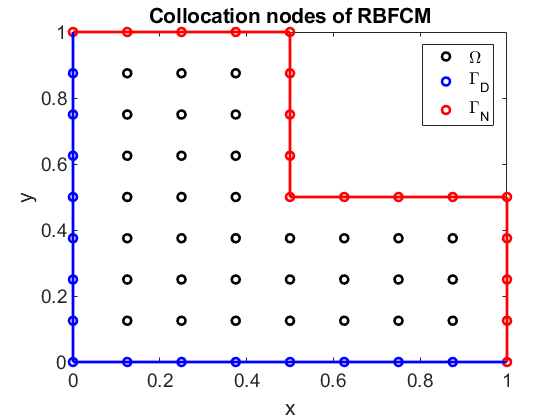}
  \end{minipage}
    \begin{minipage}[t]{0.47\textwidth}
    \vspace{0pt}
    \includegraphics[width=1.0\textwidth]{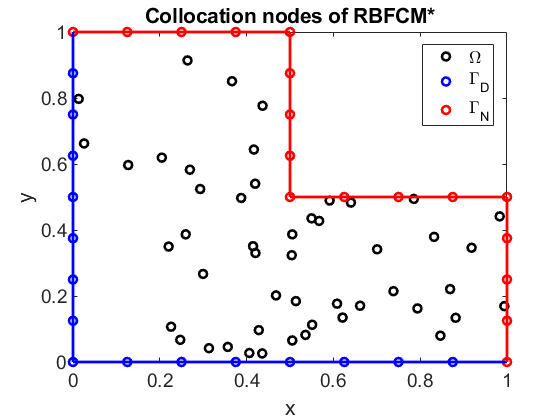}
  \end{minipage}
\end{center}
\caption{Computational domains for Example (\ref{lshapedss}): FEM O(1) elements (upper left), FEM O(2) elements (upper right), uniform RBFCM nodes (lower left), randomly generated RBFCM nodes (lower right) $\delta h=0.125$.}
  \label{fig:lshapedssdomain-set-up}
\end{figure}

In this example, the domain is changed to an L-shaped domain, $\Omega=(0,1)\times(0,1)\setminus[0.5,1)\times[0.5,1) \subset \mathbb{R}^2$, and the boundary conditions are of the Dirichlet, $\Gamma_D=\{(0,y)\cup(x,0)\}\subseteq\partial\Omega$, and the Neumann type, $\Gamma_N= \partial\Omega \setminus \Gamma_D$. 
In order to show the true meshless nature of RBFCM, an additional domain, in which the nodes are randomly distributed, is created along with the regular uniform domains (see Figure \ref{fig:lshapedssdomain-set-up}). The number of nodes used in this non-uniform domain is the same as its uniform counterpart. 
The results obtained by using this non-uniform domain are marked with an asterisk sign \textit{(e.g., RBFCM-MQ*)}. The analytical solution is set to be:
\begin{align}
\label{eq:lshapedssanalytic}
u(x,y)=ye^{-x^2}+ \sin(\pi x) \cos(\pi y)
\end{align}
Thus, by setting $k=1$ in in Equation \ref{eq:Poisson_strongnumformfem}, the source term, $f$ becomes:
\begin{align}
f=&ye^{-x^2}(2-4x^2)+2\pi^2\sin(\pi x) \cos(\pi y)
\end{align}
The analytical solution used in this example is somewhat similar to that of Example \ref{example1}. However, it is significant to see the effect of the Neumann boundary condition on the approximation behaviors of the methods.

\begin{table}[H]
\caption{Domain specifications of Example \ref{lshapedss}.}
\centering
\begin{tabular}{c m{4cm} m{4cm} m{4cm}}
\hline
\boldmath$\delta h$&\textbf{FEM O(1)}&\textbf{FEM O(2)}&\textbf{RBFCM(MQ and TPS)}\\
\hline
1/4&5 nodes in $\Omega$, 7 nodes on $\Gamma_D$, 9 nodes on $\Gamma_N$, 24 elements.&33 nodes in $\Omega$, 15 nodes on $\Gamma_D$, 17 nodes on $\Gamma_N$, 24 elements.&5 nodes in $\Omega$, 7 nodes on $\Gamma_D$, 9 nodes on $\Gamma_N$.\\
1/8&33 nodes in $\Omega$,  15 nodes on $\Gamma_D$, 17 nodes on $\Gamma_N$, 96 elements.&161 nodes in $\Omega$, 31 nodes on $\Gamma_D$, 33 nodes on $\Gamma_N$, 96 elements.&33 nodes in $\Omega$, 15 nodes on $\Gamma_D$, 17 nodes on $\Gamma_N$.\\
1/16&161 nodes in $\Omega$, 31 nodes on $\Gamma_D$, 33 nodes on $\Gamma_N$, 384 elements.& 705 nodes in $\Omega$, 63 nodes on $\Gamma_D$, 65 nodes on $\Gamma_N$, 384 elements. &161 nodes in $\Omega$, 31 nodes on $\Gamma_D$, 33 nodes on $\Gamma_N$.\\
1/32&705 nodes in $\Omega$, 63 nodes on $\Gamma_D$, 65 nodes on $\Gamma_N$, 1536 elements.& 2945 nodes in $\Omega$, 125 nodes on $\Gamma_D$, 127 nodes on $\Gamma_N$, 1536 elements. &705 nodes in $\Omega$, 63 nodes on $\Gamma_D$, 65 nodes on $\Gamma_N$.\\
\hline
\end{tabular}
\label{dse32}
\end{table}

Considering the LSE results, the convergence achieved by both RBFCM-TPS and RBFCM-TPS* is very satisfactory, and RBFCM-TPS* outperforms FEM O(2) for smaller mesh sizes. The LSE results of RBFCM-MQ, on the other hand, performed less well in terms of accuracy. It was not able to provide more accurate results than FEM O(2) even for coarser meshes. This might be due to the extremely large condition number of the system matrices generated by the RBFCM-MQ. The large condition number escalated the inversion errors and worsened the approximation. Also, no improvement is observed in RBFCM-MQ* approximation.

\begin{table}[H] 
\caption[Comparison results for Example \ref{lshapedss}.]{Comparison results for Example \ref{lshapedss}. (LSE: Least Square Error; RMSE: Root Mean Square Error; MRE: Maximum Relative Error; CN: Condition Number; RT: Runtime (sec); OSP: Optimum Shape Parameter. Asterisk sign, "*", indicates RBFCM solutions where the nodes are randomly distributed.)}
\begin{center}
\begin{tabular}{|c|c|c|c|c|c|c|c|} \hline
	\boldmath$\delta h$&& \textbf{LSE} & \textbf{RMSE}& \textbf{MRE}&\textbf{CN}& \textbf{RT}&\textbf{OSP} \\
\hline
	1/4&FEM O(1) & 6.131e-02 & 5.268e-02 &4.955e-01& 8.547e+16 & 0.05&\\
	 &FEM O(2) & 3.973e-03 & 2.622e-03 & 4.511e-02& 2.271e+16 & 4.20&\\
	 &RBFCM-MQ & 1.048e-02 & 1.160e-02 & 2.573e-01& 9.793e+38 & 1.31&2.045\\
	 &RBFCM-TPS & 3.961e-01 & 4.307e-01 & 6.552e-00& 1.382e+04 & 0.02&\\
	 &RBFCM-MQ* & 1.270e-02 & 1.834e-02 & 3.108e-01& 1.026e+15 & 1.31&3.914\\
	 &RBFCM-TPS* & 1.046e-01 & 7.277e-01 & 1.092e+01 & 1.975e+04 & 0.02&\\
	 \hline
	 	1/8&FEM O(1) & 1.681e-02 & 1.265e-02 & 6.537e-00  & 4.650e+16 &0.62 &\\
	 &FEM O(2) & 4.985e-04 & 2.139e-04 & 4.495e-02& 3.147e+16 & 16.71&\\
	 &RBFCM-MQ & 8.401e-04 & 1.049e-03 & 7.801e-03& 2.851e+76 & 13.25&4.985\\
	 &RBFCM-TPS & 5.360e-03 & 7.809e-03 & 1.717e-00& 2.888e+08 & 0.14 &\\
	 &RBFCM-MQ* & 5.477e-04 & 3.712e-04 & 6.387e-03& 1.240e+12 & 13.61&0.71\\
	 &RBFCM-TPS* & 1.074e-03 & 3.346e-02 & 6.832e-01& 1.167e+08 & 0.12&\\
	 \hline
	1/16&FEM O(1) & 4.317e-03 & 3.048e-03 & 1.622e+00& 8.969e+16 &2.55&  \\
	 &FEM O(2) & 6.262e-05 & 1.797e-05 & 1.434e-02& 2.384e+16 & 70.69&\\
	 &RBFCM-MQ & 1.254e-04 & 1.491e-04 & 2.294e-02& 5.676e+131 & 171.82&3.348\\
	 &RBFCM-TPS & 1.247e-04 & 1.639e-04 & 7.994e-02& 5.861e+11 & 1.57 &\\
	 &RBFCM-MQ* & 6.728e-04 & 2.765e-04 & 1.297e-02& 6.293e+20 & 174.13&2.819\\
	 &RBFCM-TPS* & 2.078e-05 & 1.001e-03 & 2.072e-01& 1.545e+12 & 1.59&\\
	\hline
		1/32&FEM O(1) & 1.087e-03 & 7.421e-04 & 9.782e-01 & 7.573e+16 &11.22&\\
	 &FEM O(2) & 7.856e-06 & 1.546e-06 & 5.908e-04& 1.065e+17 & 315.91&\\
	 &RBFCM-MQ & 5.946e-05 & 6.029e-05 & 5.981e-02& 5.893e+164 & 2525.76& 2.127\\
	 &RBFCM-TPS & 8.116e-06 & 1.054e-05 & 4.855e-03& 7.024e+14 & 22.02 &\\
	 &RBFCM-MQ* & 1.377e-04 & 1.422e-04 & 9.717e-02& 8.7531e+16 & 2531.83&0.182\\
	 &RBFCM-TPS* & 4.378e-06 & 5.619e-05 & 6.827e-03& 8.237e+15 & 22.09&\\
	 \hline
	\end{tabular}
\label{tab:ex2ErrComp}
\end{center}
\end{table}

\begin{figure}[H]
  \begin{center}
    \begin{minipage}[t]{0.47\textwidth}
    \vspace{0pt}
    \includegraphics[width=1.0\textwidth]{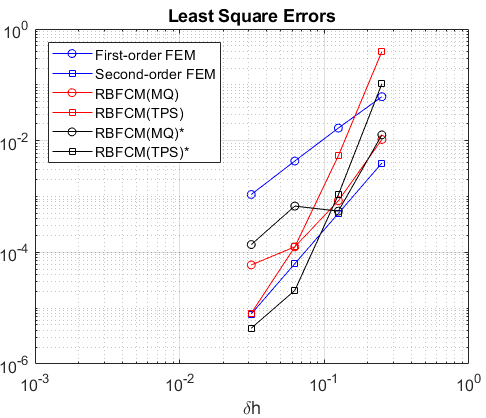}
  \end{minipage}
  \hspace{0.02\textwidth}
  \begin{minipage}[t]{0.47\textwidth}
    \vspace{0pt}
    \includegraphics[width=1.0\textwidth]{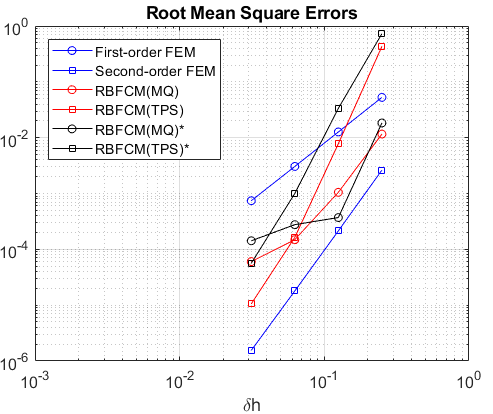}
  \end{minipage}
\end{center}
\caption{Least square errors (left), and root mean square errors (right) of Example \ref{lshapedss}.}
  \label{fig:lshapedh-conv}
\end{figure}

\subsection{Unsteady Poisson Equation with Dirichlet Boundary Conditions}
\label{example3}
The domain and the boundary conditions in this example are the same as in Example \ref{example1}.\footnote{The reader may refer to Table \ref{dse31} for the domain specifications.} Time step size, $\delta t$, is set to be 0.01. All the results and figures related to this example are of the final time step, $T_f=50$. 
An analytical solution, in which the time term is linear, is defined here to eliminate any error related to the time integrator since Backward Euler Method we are using to discretize the time term is only first-order accurate. The analytical solution is then written as follows:
\begin{align}
\label{eq:ex3analytic}
u(x,y,t)= \exp{\left( \frac{-x}{y+1} \right)}+0.8t
\end{align}
Therefore, by setting $k=1$, the source term $f$ becomes:
\begin{align}
f=0.8-\frac{\exp{{\left( \frac{-x}{y+1} \right)}}}{(y+1)^2}-\frac{x\exp{{\left( \frac{-x}{y+1} \right)}(x-2y-2)}}{(y+1)^4}\notag
\end{align}

\begin{table}[H]
\caption[Comparison results for Example \ref{example3}.]{Comparison results for Example \ref{example3}. (LSE: Least Square Error; RMSE: Root Mean Square Error; MRE: Maximum Relative Error; CN: Condition Number; RT: Runtime (sec); OSP: Optimum Shape Parameter.)}
\begin{center}
\begin{tabular}{|c|c|c|c|c|c|c|c|} \hline
	\boldmath$\delta h$&& \textbf{LSE} & \textbf{RMSE}& \textbf{MRE}&\textbf{CN}& \textbf{RT}&\textbf{OSP} \\
	\hline
	1/4&FEM O(1) & 2.981e-03 & 3.438e-04 &6.716e-04& 1.472e+00 & 0.26&\\
	 &FEM O(2) & 8.227e-05 & 2.367e-05 & 6.7724e-05& 6.447e+00 & 7.74&\\
	 &RBFCM-MQ & 2.243e-04 & 1.245e-04 & 3.405e-04& 9.269e+13 & 3.20&3.108\\
	 &RBFCM-TPS & 5.700e-02 & 3.152e-02 & 5.830e-02& 5.683e+02 & 0.03 &\\
	 \hline
	 	1/8&FEM O(1) & 7.373e-04 & 1.018e-04 & 1.899e-04 & 4.596e+00 &1.14&\\
	 &FEM O(2) & 9.705e-06 & 1.874e-06 & 7.031e-06& 2.349e+01 & 31.34&\\
	 &RBFCM-MQ & 1.620e-04 & 6.176e-05 & 1.962e-04& 4.515e+19 & 23.72&3.790\\
	 &RBFCM-TPS & 2.094e-03 & 7.951e-04 & 1.551e-03& 2.733e+04 & 0.21 &\\
	 \hline
	1/16&FEM O(1) & 1.837e-04 & 2.731e-05 & 4.873e-05& 1.741e+01 &6.11& \\
	 &FEM O(2) & 1.186e-06 & 1.304e-07 & 6.028e-07& 9.190e+01 & 156.61&\\
	 &RBFCM-MQ & 1.555e-05 & 1.267e-05 & 4.009e-05& 1.010e+12 & 257.45&0.362\\
	 &RBFCM-TPS & 1.507e-04 & 4.990e-05 & 2.114e-04& 1.767e+06 & 2.72 &\\
	\hline
		1/32&FEM O(1) & 4.587e-05 & 7.059e-06 & 1.222e-05 & 6.872e+01 &51.58&\\
	 &FEM O(2) & 1.472e-07 & 8.557e-09 & 4.458e-08& 3.656e+02 & 1133.21&\\
	 &RBFCM-MQ & 7.028e-06 & 6.423e-06 & 2.519e-05& 1.584e+12 & 3510.49&0.168\\
	 &RBFCM-TPS & 1.428e-05 & 7.325e-06 & 6.969e-05& 1.148e+08 & 40.04 &\\
	 \hline
	\end{tabular}
\label{tab:ex3ErrComp}
\end{center}
\end{table}

\begin{figure}[H]
  \begin{center}
    \begin{minipage}[t]{0.47\textwidth}
    \vspace{0pt}
    \includegraphics[width=1.0\textwidth]{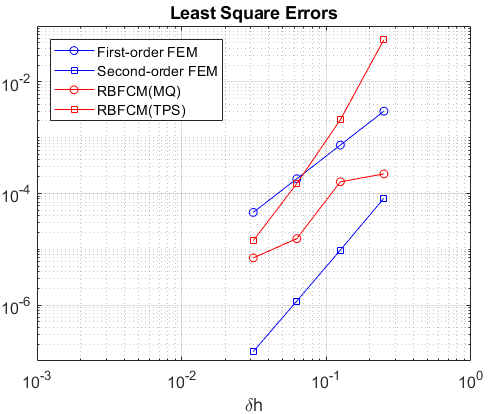}
  \end{minipage}
  \hspace{0.02\textwidth}
  \begin{minipage}[t]{0.47\textwidth}
    \vspace{0pt}
    \includegraphics[width=1.0\textwidth]{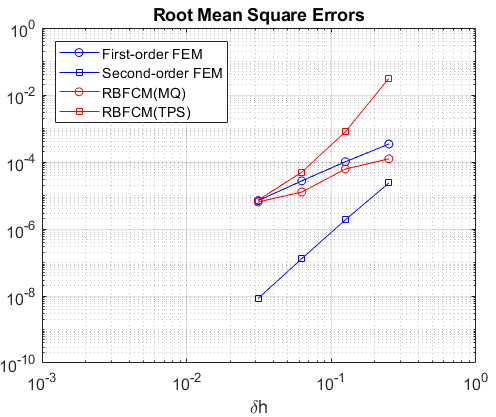}
  \end{minipage}
\end{center}
\caption{Least square error (left), and root mean square error (right) of Example \ref{example3} at the final time step.}
  \label{fig:ex3h-conv}
\end{figure}

Looking at Table \ref{tab:ex3ErrComp} and Figure \ref{fig:ex3h-conv} and considering both LSE and RMSE, the accuracy of FEM O(2) is much better than any other alternative, especially for finer meshes. 
However, as in the previous examples, the convergence rate of RBFCM-TPS is noteworthy even though its accuracy for the case that the least number of nodes used relatively low.

\subsection{Steady Stokes Equation with Dirichlet Boundary Conditions - Colliding Flow}
\label{colliding}

In this colliding flow example, we are dealing with a rather complex PDE, the Stokes equation. The domain is $\Omega=(-1,1)\times(-1,1)\subset \mathbb{R}^2$. Rather than using natural boundary conditions at the outflow, which is a natural choice for Stokes equation, the Dirichlet boundary conditions are defined for the velocity at all boundaries ($\Gamma_D=\partial \Omega$).\footnote{Even though the domain of this example is bigger than that of Example \ref{example1}, the reader may refer to Table \ref{dse31} and think of the values in $\delta h$ column as their double (\textit{i.e.}, $1/4$ as $1/2$).} 
This leads to a situation where pressure is determined only up to a constant. 
Hence, the pressure is not uniquely defined (see \cite{wathenbook} for details). 
Because of this issue, special attention must be given to pressure approximation. 
A couple of techniques to overcome this problem is available in the literature (see \cite{brezzifortin}). Even though it has no physical meaning,  rather than satisfying continuity at the upper right node($x,y$=(1,1)), the exact value of the pressure is forced as a Dirichlet boundary condition in order to obtain a unique numerical solution for pressure unknown. This technique is used for both FEM and RBFCM. As a result, the pressure field is uniquely determined where it is compatible with the analytical solution, rather than shifting up or down around the analytical solution.
Analytical solutions for the unknowns of colliding flow are:
\begin{align}
	u_x=20xy^3, \textnormal{ }u_y=5x^4-5y^4,\textnormal{ } p=60x^2y-20y^3
\end{align}
LSE results show that RBFCM-MQ starts with higher accuracy for both velocity and pressure. However, as we refine the mesh, the velocity approximation of FEM outperforms RBFCM-MQ because of its higher convergence rate. RBFCM MQ, according to our results, shows no convergence pattern for solving colliding flow, and its condition number problem continues. To approximate pressure, FEM O(1) is used. Therefore, while optimizing the shape parameter for RBFCM-MQ, velocity errors are prioritized to have comparable results between the two methods.

\begin{figure}[H]
  \begin{center}
    \begin{minipage}[t]{0.47\textwidth}
    \vspace{0pt}
    \includegraphics[width=1.0\textwidth]{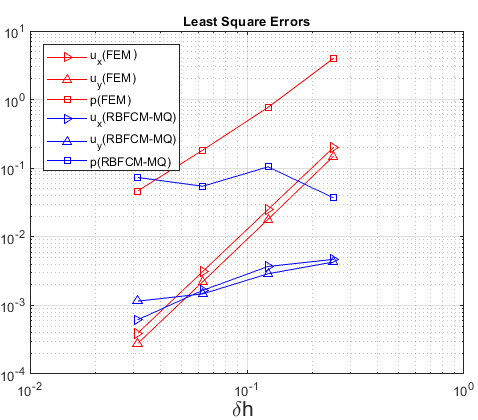}
  \end{minipage}
  \hspace{0.02\textwidth}
  \begin{minipage}[t]{0.47\textwidth}
    \vspace{0pt}
    \includegraphics[width=1.0\textwidth]{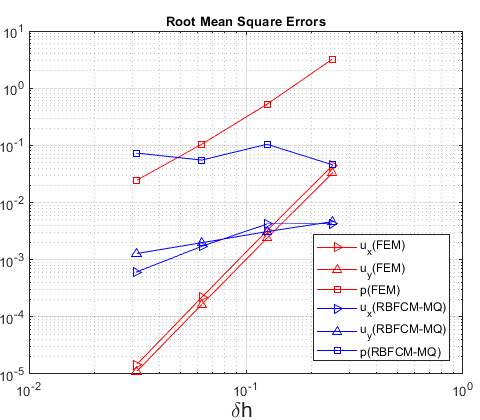}
  \end{minipage}
\end{center}
\caption{Least square errors (left), and root mean square errors (right) of Example \ref{colliding}.}
  \label{fig:ex2StokesErr}
\end{figure}

\begin{table}[H] 
\caption[Comparison results for Example \ref{colliding}.]{Comparison results for Example \ref{colliding}. (FEM: FEM O(2) for velocity and FEM O(1) for pressure approximations; LSE: Least Square Error; RMSE: Root Mean Square Error; MRE: Maximum Relative Error; CN: Condition Number; RT: Runtime (sec); OSP: Optimum Shape Parameter. Asterisk sign, "*", indicates RBFCM solutions where the nodes are randomly distributed.)}
\begin{center}
\begin{tabular}{|c|c|c|c|c|c|c|c|c|} \hline
\boldmath$\delta h$&&& \textbf{LSE} & \textbf{RMSE}& \textbf{MRE}&\textbf{CN}& \textbf{RT}&\textbf{OSP} \\
	\hline
	1/4&FEM &$u_x$& 1.997e-01 & 4.431e-02 & 6.274e-01&&&\\
	 & &$u_y$& 1.997e-01 & 3.258e-02 & 1.327e-01&  2.903e+04 & 6.63&\\
	 & &$p$& 3.381e-00 & 3.225e-00 & 4.9761e-01&&&\\
	 \cline{2-9}
	 &RBFCM-MQ &$u_x$& 4.694e-03 & 4.230e-03 & 6.104e-03&&&\\
	 & &$u_y$& 4.290e-03 & 4.618e-03 & 1.650e-02&  1.158e+20 & 1.15 & 34.978\\
	 & &$p$& 3.702e-02 & 4.572e-02 & 1.028e-02&&&\\
	 \hline
	 \hline
	 	1/8&FEM &$u_x$& 2.499e-02 & 3.214e-03 & 5.257e-01&&&\\
	 & &$u_y$& 1.796e-02 & 2.381e-03 & 1.170e-01&  3.919e+05 &27.42&\\
	 & &$p$& 7.713e-01 & 5.242e-01 & 3.803e-01&&&\\
	 \cline{2-9}
	 	 &RBFCM-MQ &$u_x$& 3.688e-03 & 4.209e-03 & 4.375e-02&&&\\
	 & &$u_y$& 2.892e-03 & 3.103e-03 & 3.383e-01&  1.502e+21 & 12.08 & 31.781\\
	 & &$p$& 1.046e-01 & 1.043e-01 & 3.470e-01&&&\\
	 \hline
	 \hline
	 	1/16&FEM &$u_x$& 3.118e-03 & 2.189e-04 & 5.247e-01&&&\\
	 & &$u_y$& 2.215e-03 & 1.637e-04 & 1.168e-01& 5.672e+06 &126.63 &\\
	 & &$p$& 1.833e-01 & 1.046e-01 &3.349e-01&&&\\
	 \cline{2-9}
	 	 &RBFCM-MQ &$u_x$& 1.659e-03 & 1.710e-03 & 3.078e-01&&&\\
	 & &$u_y$& 1.471e-03 & 1.980e-03 & 1.000e-00&  2.101e+22 & 183.94 & 7.516\\
	 & &$p$& 5.458e-02 & 5.526e-02 &1.365e-00&&&\\
	 \hline
	 \hline
	 	1/32&FEM &$u_x$& 3.893e-04 & 1.445e-05 & 5.248e-01&&&\\
	 & &$u_y$& 2.756e-04 & 1.094e-05 & 1.173e-01& 8.585e+07 &1587.05&\\
	 & &$p$& 4.627e-02 & 2.431e-02 & 9.817e-01&&&\\
	 \cline{2-9}
	 	&RBFCM-MQ &$u_x$& 6.143e-04 & 6.069e-04 & 1.442e-00&&&\\
	 & &$u_y$& 1.152e-03 & 1.277e-03 & 3.556e-00&  7.188e+21 & 4549.53 & 2.034\\
	 & &$p$& 7.330e-02 & 7.335e-02 & 1.513e+01&&&\\
	 \hline
	\end{tabular}
\label{tab:exS1Err}
\end{center}
\end{table}
\pagebreak
\subsection{Unsteady Stokes Equation with Natural Boundary Conditions}
\label{USEwNBC}
A time-dependent Stokes Equation is modeled in this example. The domain is chosen to be L-shaped, and it is exactly identical to that of Example \ref{lshapedss}. The RBFCM solutions for randomized nodes are also present here. Time step size, $\delta t$, is chosen to be 0.01, and all the results are obtained at  the final time step, $T_f=50$. The natural boundary condition is applied at the right boundary, $\Gamma_N=(1,y)\subset\partial\Omega$, and the rest of the boundaries are the Dirichlet type, $\Gamma_D=\partial\Omega\setminus\Gamma_N$. This example has very basic analytical solutions. The only time-dependent variable here is the velocity in the x-direction, $u_x$:
\begin{align}
	u_x=t+1-y^3, \textnormal{ }u_y=-x^3+3x^2-3x,\textnormal{ } p=-6xy-x+6y+constant
\end{align}
Here, the constant term in the pressure is set to $1$ to apply natural boundary conditions at the boundary $x=1$. By doing this, the boundary term on the right hand side of Equation \ref{weakStokes} becomes $0$ at $x=1$, meaning that at the natural boundary, the velocity flux and pressure are at balance.

\begin{figure}[H]
  \begin{center}
    \begin{minipage}[t]{0.47\textwidth}
    \vspace{0pt}
    \includegraphics[width=1.0\textwidth]{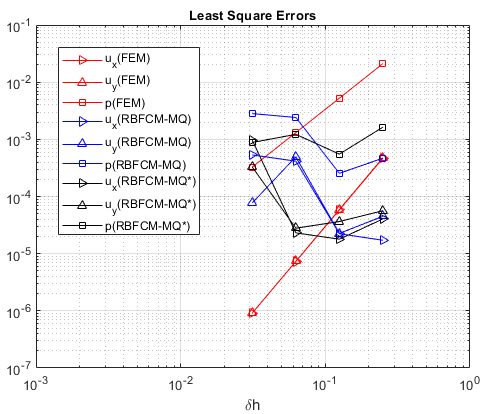}
  \end{minipage}
  \hspace{0.02\textwidth}
  \begin{minipage}[t]{0.47\textwidth}
    \vspace{0pt}
    \includegraphics[width=1.0\textwidth]{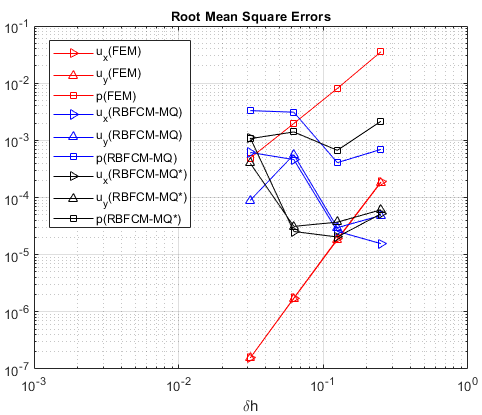}
  \end{minipage}
\end{center}
\caption{Least square errors (left), and root mean square errors (right) of Example \ref{USEwNBC} at the final time step, $T_f=50$.}
  \label{fig:unStokesEx1Errs}
\end{figure}

This example additionally includes an RBFCM-MQ model with non-uniform, randomized nodes in the domain, RBFCM-MQ*. However, there is not a significant improvement in the error results. The approximation of FEM O(2) was better in the previous unsteady Poisson problem, and it converged better for all mesh sizes. In this example, RBFCM-MQ and RBFCM-MQ* approximate the analytical solution better for the biggest $\delta h$ value (\textit{i.e.}, for coarsest mesh). Anyhow, an excellent match with expected convergence behavior is observed for FEM, whereas RBFCM did not show any convergence pattern, and besides, its accuracy becomes worse as the mesh gets finer.

\begin{table}[H]
\caption[Comparison results for Example \ref{USEwNBC}.]{Comparison results for Example \ref{USEwNBC}. (FEM: FEM O(2) for velocity and FEM O(1) for pressure approximations; LSE: Least Square Error; RMSE: Root Mean Square Error; MRE: Maximum Relative Error; CN: Condition Number; RT: Runtime (sec); OSP: Optimum Shape Parameter.)}
\begin{center}
\begin{tabular}{|c|c|c|c|c|c|c|c|c|} \hline
\boldmath$\delta h$&&& \textbf{LSE} & \textbf{RMSE}& \textbf{MRE}&\textbf{CN}& \textbf{RT}&\textbf{OSP} \\
	\hline
	 1/4&FEM &$u_x$& 4.663e-04 & 1.876e-04 & 2.396e-04&&&\\
	 & &$u_y$& 4.767e-04 & 1.833e-04 & 1.2356e-03&  4.235e+04 & 6.93&\\
	 & &$p$& 2.115e-02 & 3.539e-02 & 4.976e-01&&&\\
	 \cline{2-9}
	 &RBFCM-MQ &$u_x$& 1.723e-05 & 1.560e-05 & 2.216e-05&&&\\
	 & &$u_y$& 4.509e-05 & 4.842e-05 & 2.216e-05&  5.039e+21 & 1.26 & 29.682\\
	 & &$p$& 4.680e-04 & 6.937e-04 & 3.204e-03&&&\\
	 \cline{2-9}
	 &RBFCM-MQ* &$u_x$& 4.059e-05 & 5.171e-05 & 5.929e-05&&&\\
	 & &$u_y$& 5.637e-05 & 6.087e-05 & 1.761e-04&  6.262e+20 & 1.28 & 31.906\\
	 & &$p$& 1.632e-03 & 2.154e-03 & 3.339e-03&&&\\
	 \hline
	 \hline
	 1/8&FEM &$u_x$& 5.817e-05 & 1.885e-05 & 3.180e-05&&&\\
	 & &$u_y$& 5.912e-05 & 1.821e-05 & 1.621e-04&  1.292e+05 &31.18&\\
	 & &$p$& 5.239e-03 & 8.089e-03 & 7.590e-02&&&\\
	 \cline{2-9}
	 &RBFCM-MQ &$u_x$& 2.232e-05 & 2.559e-05 & 2.969e-05&&&\\
	 & &$u_y$& 2.264e-05 & 2.969e-05 & 8.234e-05&  6.391e+20 & 9.58 & 33.947\\
	 & &$p$& 2.544e-04 & 4.068e-04 & 1.198e-03&&&\\
	 \cline{2-9}
	 &RBFCM-MQ* &$u_x$& 1.790e-05 & 2.039e-05 & 4.619e-05&&&\\
	 & &$u_y$& 2.174e-05 & 3.704e-05 & 9.483e-04&  3.586e+20 & 9.46 & 33.640\\
	 & &$p$& 5.513e-04 & 6.690e-04 & 1.9989e-02&&&\\
	 \hline
	 \hline
	 1/16&FEM &$u_x$& 7.278e-06 & 1.741e-06 & 3.516e-06&&&\\
	 & &$u_y$& 7.347e-06 & 1.713e-06 & 2.093e-05& 4.414e+05 &216.46 &\\
	 & &$p$& 1.310e-03 & 1.971e-03 &3.717e-02&&&\\
	 \cline{2-9}
	 &RBFCM-MQ &$u_x$& 4.146e-04 & 4.612e-04 & 8.323e-04&&&\\
	 & &$u_y$& 4.930e-04 & 5.627e-04 & 5.307e-03&  2.204e+22 & 107.74 & 43.600\\
	 & &$p$& 2.411e-03 & 3.104e-03 &7.526e-02&&&\\
	 \cline{2-9}
	 &RBFCM-MQ* &$u_x$& 2.311e-05 & 2.544e-05 & 7.281e-05&&&\\
	 & &$u_y$& 1.054e-04 & 3.126e-05 & 9.763e-04&  1.479e+22 & 109.44 & 32.305\\
	 & &$p$& 1.223e-03 & 1.402e-03 & 4.033e-02&&&\\
	 \hline
	 \hline
	 1/32&FEM &$u_x$& 9.106e-07 & 1.564e-07 & 3.953e-07&&&\\
	 & &$u_y$& 9.152e-07 & 1.559e-07 & 2.638e-06& 1.623e+06 &3088.72&\\
	 & &$p$& 3.275e-04 & 4.895e-04 & 1.849e-02&&&\\
	 \cline{2-9}
	 &RBFCM-MQ &$u_x$& 5.336e-04 & 6.132e-04 & 9.326e-04&&&\\
	 & &$u_y$& 7.722e-05 & 8.919e-05 & 1.170e-03&  3.983e+22 & 2314.97 & 6.091\\
	 & &$p$& 4.003e-03 & 4.754e-03 & 2.011e-02&&&\\
	 \cline{2-9}
	 &RBFCM-MQ* &$u_x$& 9.638e-04 & 1.105e-03 & 2.354e-03&&&\\
	 & &$u_y$& 3.284e-04 & 3.963e-04 & 2.762e-02&  1.165e+24 & 2317.63 & 6.841\\
	 & &$p$& 8.797e-04 & 1.072e-03 & 5.911e-01&&&\\
	 \hline
	\end{tabular}
\label{tab:exS2Err}
\end{center}
\end{table}

\section{Discussions}

The FEM and RBCFM performances regarding the accuracy, runtime, stability, and the ease of implementation criteria for solving the steady and unsteady Poisson and Stokes equations were investigated to further the basis for comparing the two methods. The accuracy criteria were based mostly on the least square error over numerical examples where the analytical solutions were known a priori.

Runtime differs with the model of the central processing unit (CPU) used. The numerical models in this work are processed with Intel Core i5-7200U @ 2.50 GHz. The software that the models are executed is Python v3.7. System matrix inversions of all models are calculated by "pinv" command of numpy module. This command computes Moore-Penrose inversion using singular-value decomposition of the matrix. The preprocessing and the postprocessing time are not included in runtime.

The following observations and conclusions were made from the numerical experiments:
\begin{itemize}
	\item In Poisson cases:
	\begin{enumerate}[label=$\circ$]
		\item RBFCM-MQ solutions yielded more accurate results with respect to FEM O(1) solutions.
		\item RBFCM-TPS solutions started with higher errors then reached the accuracy level of FEM O(1) solutions in all of the examples with decreasing mesh size.
		\item In general, FEM O(2) solutions provided more accurate solutions for all the cases irrespective of the mesh size because of the fact that the FEM O(2) utilizes more nodes for the same mesh size (see Tables \ref{dse31},\ref{dse32}).
		\item Condition numbers of the system matrices generated by RBFCM in general, by RBFCM-MQ in particular, are high.
		\item System matrices of RBFCM-TPS are better conditioned than RBFCM-MQ. On the other hand, both FEM O(1) and O(2) system matrices have proper condition numbers as they should be.
	\end{enumerate}
	\item In Stokes cases:
	\begin{enumerate}[label=$\circ$]
		\item RBFCM-MQ solutions are more accurate for coarser meshes (\textit{i.e.}, bigger $\delta h$ values) compared to FEM solutions.
		\item Refining the mesh does not necessarily improve the accuracy of RBFCM-MQ approximations.
		\item The velocity approximations obtained show that FEM follows an expected convergence pattern, whereas no meaningful pattern in RBFCM-MQ approximation was observed, and as a result, FEM O(2) becomes more accurate as the mesh is refined. 
	\end{enumerate}
	\item Error convergence behavior with respect to mesh size assessment shows that:
	\begin{enumerate}[label=$\circ$]
		\item Both first and FEM O(2) solution errors linearly decrease on a log-log scale. The slope of the LSE is optimal, and as expected, the slope is $2$ for the  FEM O(1) solutions and $3$ for FEM O(2) solutions.
		\item Even though the convergence behavior of RBFCM-TPS models is not well defined as FEM models, it was satisfactory for Poisson examples.
		\item The error convergence behavior of RBFCM-MQ models is poor and unpredictable since using more nodes in our models did not necessarily improve the accuracy. 
	\end{enumerate}
	\item Runtime assessment shows that:
	\begin{enumerate}[label=$\circ$]
		\item Most of the runtime (around 98-99\%) spent in RBFCM-MQ models is due to the shape parameter optimization algorithm because the algorithm seeks the optimum shape parameter by running the model several times using different shape parameters and comparing the RMSE for each shape parameter iteratively. It is possible to use shape parameters proposed in \cite{Hardymq, franke1975, konghonetal, wuhon2003} without seeking an optimum one.
		\item If the literature-available shape parameters are used without the iterative scheme, the runtime of RBFCM-MQ would be very close to that of RBFCM-TPS.
		\item For coarser meshes (\textit{i.e.}, higher $\delta h$), RBFCM-TPS models run faster than both first and FEM O(2) solutions. As the mesh is refined, however, FEM O(1) approximation becomes faster. 
		\item For the same mesh size, FEM O(2) provides more information to the solver (see Tables \ref{dse31}, \ref{dse32}). Thus, the time needed to execute FEM O(2) models is longer. 
	\end{enumerate}
	\item Ease of use shows that:
	\begin{enumerate}[label=$\circ$]
		\item The implementation of RBFCM is far simpler. The model only needs the derivatives of radial basis functions and the locations of the nodes distributed in the domain.
		\item To model FEM, on the other hand, a connectivity matrix, an external integrator such as Gaussian quadrature, element transformations, calculation of Jacobians, etc. are needed. In the modeling of L-shaped examples where the nodes are randomly distributed, no additional labor is needed for the RBFCM model. However, creating a mesh for the same domain with FEM would not be as straightforward.
	\end{enumerate}
\end{itemize}

\begin{figure}[H]
  \begin{center}
    \begin{minipage}[t]{0.30\textwidth}
    \vspace{0pt}
    \includegraphics[width=1.0\textwidth]{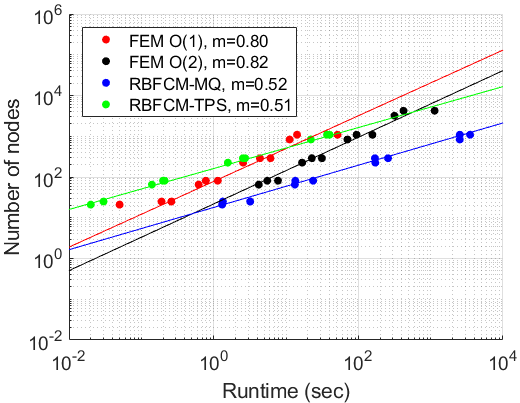}
  \end{minipage}
  \hspace{0.02\textwidth}
  \begin{minipage}[t]{0.30\textwidth}
    \vspace{0pt}
    \includegraphics[width=1.0\textwidth]{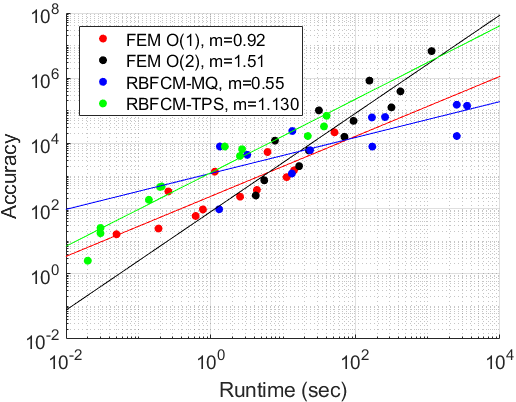}
  \end{minipage}
  \hspace{0.02\textwidth}
  \begin{minipage}[t]{0.30\textwidth}
    \vspace{0pt}
    \includegraphics[width=1.0\textwidth]{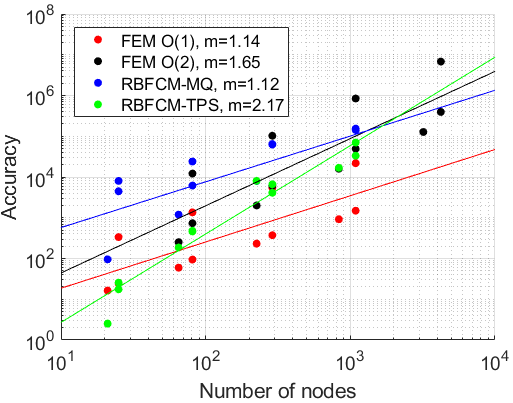}
  \end{minipage}
\end{center}
\caption{Trends for Poisson problem. Number of nodes vs. runtime at the left, accuracy vs. runtime at the middle, accuracy vs. number of nodes at the right. (Accuracy=LSE$^{-1}$; $m$ is the slope of the trend lines)}
  \label{fig:discussionPois}
\end{figure}

To have a wider overview, all FEM and RBFCM approximations are put together, and the trends are observed for number of nodes vs. runtime, accuracy vs. runtime, and accuracy vs. number of nodes. These plots provide guidance on comparing the two numerical approaches when solution requirements are set a priori. The trends noted for the Poisson equation solutions are shown in Figure \ref{fig:discussionPois}.

\begin{itemize}
	\item Both first and FEM O(2) requires less additional runtime for the same amount of increase in the number of nodes, meaning that it is less costly, in terms of runtime, for FEM solutions to increase the number of nodes.
	\item Depending on the level of accuracy required, the preferred method may require less runtime (RBFCM-MQ for less accuracy) or more runtime (FEM O(2) for more accuracy). The accuracy trend shows that RBFCM-TPS is the most accurate method if the time is constrained up to around $10^3$ seconds. From that point on, FEM O(2)  yields the most accurate results.
	\item By refining the mesh (\textit{i.e.}, increasing the number of nodes), the accuracy of RBFCM-TPS is improved drastically. On the other hand, RBFCM-MQ provides the best accuracy for coarse meshes.
	\item The combination of the two graphs can also provide additional insight into the method to choose from. For example, if one is constricted with the number of nodes (memory) at $10^3$ while looking at running the simulations at an accuracy of $10^5$, it seems that RBFCM-TPS is the best option.
\end{itemize}

\begin{figure}[H]
  \begin{center}
    \begin{minipage}[t]{0.30\textwidth}
    \vspace{0pt}
    \includegraphics[width=1.0\textwidth]{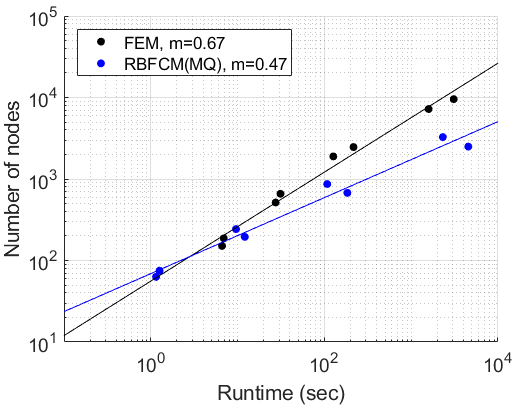}
  \end{minipage}
  \hspace{0.02\textwidth}
  \begin{minipage}[t]{0.30\textwidth}
    \vspace{0pt}
    \includegraphics[width=1.0\textwidth]{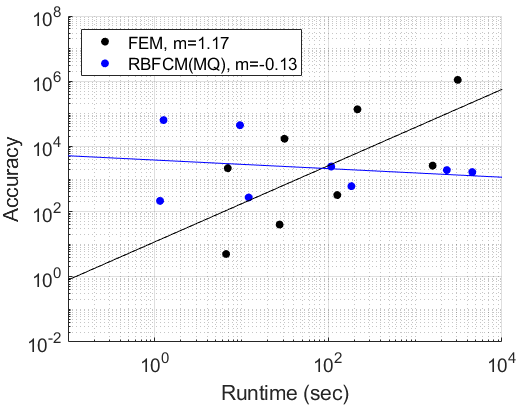}
  \end{minipage}
  \hspace{0.02\textwidth}
  \begin{minipage}[t]{0.30\textwidth}
    \vspace{0pt}
    \includegraphics[width=1.0\textwidth]{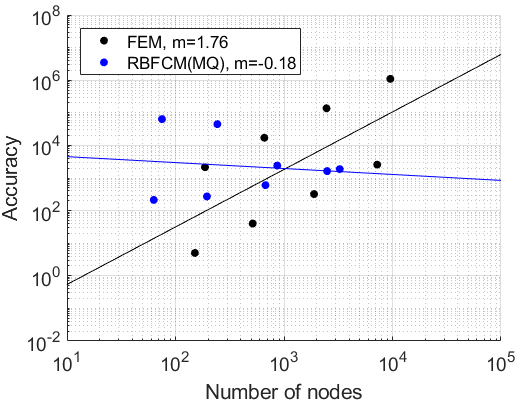}
  \end{minipage}
\end{center}
\caption{Trends for Stokes problem. Number of nodes vs. runtime at the left, accuracy vs. runtime at the middle, accuracy vs. number of nodes at the right. (Accuracy=LSE$_{u_x}^{-1}$; $m$ is the slope of the trend lines)}
  \label{fig:discussionStokes}
\end{figure}

The trends noted for the Stokes equation solutions are shown in Figure \ref{fig:discussionStokes}.

\begin{itemize}
	\item As in Poisson cases, the same amount of increase in the number of nodes cost more in RBFCM-MQ.
	\item RBFCM-MQ achieves more accurate results for less runtime or fewer number of nodes. However, the trend of accuracy is inclined downwards as the number of nodes increases. Even though it depends on the type of the problem, RBFCM-MQ does not enjoy dense node placement as reported also by several authors.
\end{itemize}

\section{Conclusions}

Comparison of numerical solutions of RBFCM and FEM for Poisson and Stokes problems considering accuracy, runtime, condition number, and ease of implementation was presented in this work. Several studies have been carried on for both FEM and RBFCM hitherto; however, comparison studies between these two methods have been limited in extent. The present research comprehensively includes comparisons over the accuracy, condition numbers, and runtimes. This work also contains observations on mesh refinement effects (\textit{i.e.}, error convergence behaviors), unsteady models, scenarios with randomly distributed RBFCM nodes, shape parameter optimization for RBFCM-MQ and its effect on runtime, models  utilizing TPS as radial basis functions along with MQ, second-order finite element approximations along with first-order finite element approximations.

Previous works have stated that RBFCM has the potential to become a solid alternative to FEM for industrial applications\cite{li2003}. The present study has identified that RBFCM models indeed can be very useful when the desired outcome is fast-running, straightforward to implement, and fairly accurate numerical model where one can grasp the solution behavior. On the other hand, FEM becomes more convenient if the desired outcome includes accuracy increase for lower cost or robust convergence behavior. However, the boundaries of these areas differ from one type of partial differential equation to another, and it should be identified through numerical tests beforehand to determine which method will satisfy the requirements needed for that case.

As a further study, authors think that extending the research to solve Navier-Stokes equations can be very interesting, as well as comparing Eulerian versus Lagrangian RBFCM models thereamong. Such a study may reveal the true potential of RBFCM since implementing the algorithm will be much more straightforward and will require less labor compared to its FEM counterpart.

\section*{Acknowledgments}
The authors gratefully acknowledge Kadir Has University for letting the resources of the university to be used throughout this study. 

\bibliographystyle{ieeetr}
\bibliography{bibby}
\end{document}